\documentclass{article}

\usepackage{amsthm}
\usepackage{amsmath}
\usepackage{amssymb}
\usepackage{mathrsfs}
\usepackage{adjustbox}
\usepackage{dsfont}
\usepackage[all]{xy}
\usepackage{tikz-cd}
\usepackage{enumitem}
\usepackage{stmaryrd}
\usepackage{hyperref}

\theoremstyle{definition}
\newtheorem{Definition}{Definition}[subsection]

\theoremstyle{plain}

\theoremstyle{plain}

\theoremstyle{plain}
\newtheorem{Proposition}[Definition]{Proposition}

\theoremstyle{plain}

\theoremstyle{plain}
\newtheorem{Lemma}[Definition]{Lemma}

\theoremstyle{plain}

\theoremstyle{plain}
\newtheorem{Corollary}[Definition]{Corollary}

\theoremstyle{plain}

\theoremstyle{plain}
\newtheorem{Conjecture}[Definition]{Conjecture}

\theoremstyle{plain}

\theoremstyle{definition}
\newtheorem{Example}[Definition]{Example}

\theoremstyle{definition}
\newtheorem{Notation}[Definition]{Notation}

\theoremstyle{remark}
\newtheorem{Remark}[Definition]{Remark}

\author{Thibault D. Décoppet}
\title{Weak Fusion 2-Categories}

\begin{document}

\bibliographystyle{alpha}




\maketitle
\hspace{1cm}
\begin{abstract}
We introduce a weakening of the notion of fusion 2-category given in \cite{DR}. Then, we establish a number of properties of (multi)fusion \linebreak 2-categories. In particular, we prove that the left and the right duals of an object in a multifusion 2-category coincide. Finally, we describe the fusion rule of the fusion 2-categories associated to certain pointed braided fusion categories.
\end{abstract}

\let\clearpage\relax

\section*{Introduction}
We present an algebraic definition of fusion 2-categories that partially fills the gap between the two existing definitions (\cite{DR} and \cite{JF}). In fact, we show that our definition is essentially equivalent to that of \cite{DR}. Further, it is possible to show that the (separable) fusion 2-categories of \cite{JF} are examples of the objects we call fusion 2-categories. The proof is not difficult, but requires developing the theory of fusion 2-categories much further, so we postpone it to \cite{D6}. However, proving the converse is a much more delicate problem. Namely, it essentially amounts to proving that every fusion 2-category is a fully dualizable object in an appropriate symmetric monoidal 4-category. We hope to return to this point in following work, but, in the meantime, we present some elementary results on the structure of (multi)fusion 2-categories.

Let us briefly describe the content of the different sections of the present article. Section \ref{sec:3categories} contains results on monoidal 2-categories that we will be needed. We begin by reviewing adjoints. Providing we are working in a monoidal 2-category that has right adjoints, we show how these adjoints can be assembled to give a monoidal 2-functor. We go on to review the definitions of duals and of coherent duals. Further, given a 2-category that has right duals, we explain how they can be put together to define a 2-functor. We conjecture that this 2-functor is in fact monoidal. Finally, we explain some properties of monoidal 2-categories that have right adjoints and right duals.

Fusion 2-categories as defined in \cite{DR} are Gray monoids satisfying some properties. In section \ref{sec:fusion}, we generalize their definition by starting with a monoidal 2-category in the sense of \cite{SP}. We also give a decomposition result for multifusion 2-categories that is analogous to that of multifusion 1-categories. Then, we show that our definition of a fusion 2-category is essentially equivalent to that of Douglas-Reutter. We continue by deriving some consequences of the existence of duals in multifusion 2-categories. Most notably, we show that right and left duals agree. Next, we show that the connected component of the unit in a fusion 2-category forms a fusion sub-2-category, and we prove that connected fusion 2-categories correspond precisely to braided fusion categories. Finally, we explain how to compute the fusion rule of the connected fusion 2-category associated to certain pointed braided fusion categories, generalizing a result of \cite{ENO}. We use this to compute the fusion rules of some fusion 2-categories.

I would like to thank Christopher Douglas and David Reutter for helpful conversations related to the content of this article.

\section{Adjoints and Duals in Monoidal 2-Categories}\label{sec:3categories}

\subsection{Adjoints}

The following definition is well-known.

\begin{Definition}
Let $\mathfrak{C}$ be a 2-category and $f:A\rightarrow B$ a 1-morphism in $\mathfrak{C}$. A right-adjoint for $f$ is a 1-morphism $f^*:B\rightarrow A$ together with 2-morphisms $\epsilon:f\circ f^*\Rightarrow Id_B$ and $\eta: Id_A\Rightarrow f^*\circ f$ satisfying the snake equations
$$(\epsilon \circ f)\cdot (f\circ \eta)=Id_f,$$
$$(f^*\circ \epsilon)\cdot (\eta\circ f^*)=Id_{f^*},$$
in which we have omitted the relevant coherence 2-isomorphisms.
One defines left-adjoints dually.
\end{Definition}

\begin{Remark}
It is well-known that a right-adjoint for a 1-morphism $f$ is unique up to unique 2-isomorphism.
\end{Remark}

\begin{Notation}
Let $\mathfrak{C}$ be a 2-category. We denote by $\mathfrak{C}^{2op}$ the 2-category obtained from $\mathfrak{C}$ by reversing the direction of the 2-morphisms. Analogously, we denote by $\mathfrak{C}^{1op}$ the 2-category obtained by reversing the direction of the 1-morphisms. If $\mathfrak{C}$ is monoidal, then it is clear that so is $\mathfrak{C}^{2op}$. The 2-category $\mathfrak{C}^{1op}$ can also be given a monoidal structure, but this construction is slightly trickier, as one needs to ``invert'' 1-morphisms. Using the algebraic definition of monoidal 2-category given in \cite{SP}, this poses no problem.
\end{Notation}

\begin{Lemma}\label{lem:adjoints2functors}
Let $\mathfrak{C}$ be a 2-category with right-adjoints. There is a 2-functor $$(-)^*:\mathfrak{C}\rightarrow \mathfrak{C}^{1op;2op}$$ that is the identity on objects and sends a 1-morphism $f$ to the 1-morphism underlying a chosen right adjoint. Dually, if $\mathfrak{C}$ has left-adjoints, there is a 2-functor denoted by $^*(-)$ that sends 1-morphisms to their left adjoints. Further, if $\mathfrak{C}$ has left and right-adjoints, $^*(-)$ is a pseudo-inverse for $(-)^*$.
\end{Lemma}
\begin{proof}
The is well-known, for instance, see \cite{Ghosh}, or appendix \ref{sub:adjointsmonoidal2functor} below.
\end{proof}

If we assume in addition that $\mathfrak{C}$ is monoidal, then the 2-functor $(-)^*$ can be made monoidal.

\begin{Lemma}\label{lem:adjointsmonoidal2functor}
Let $\mathfrak{C}$ be a monoidal 2-category with right-adjoints. The 2-functor $(-)^*$ of lemma \ref{lem:adjoints2functors} can be made monoidal. If, in addition, $\mathfrak{C}$ has left-adjoints, then the monoidal 2-functors $^*(-)$ is a monoidal pseudo-inverse for $(-)^*$.
\end{Lemma}
\begin{proof}
See appendix \ref{sub:adjointsmonoidal2functor} below.
\end{proof}

\subsection{Duals}

Some care has to be taken with respect to what one calls a dual in a monoidal 2-category $\mathfrak{C}$. For us, a right dual for an object $A$ of $\mathfrak{C}$ consists of an object $A^{\sharp}$, together with two 1-morphisms $i_A:I\rightarrow A^{\sharp}\Box A$ and $e_A:A\Box A^{\sharp}\rightarrow I$ satisying the snake equations up to 2-isomorphisms. Similarly, one can give a definition of a left dual for $A$. These definitions have the advantage of being concise and easy to check, but they are not convenient to use in constructions. That is why we shall also need to consider the refinement defined in \cite{Pstr} called a coherent right dual.

\begin{Definition}\label{def:coherentdual}
Let $A$ be an object of a monoidal 2-category $\mathfrak{C}$. A coherent right dual for $A$ consists of an object $A^{\sharp}$ in $\mathfrak{C}$, 1-morphisms $i_A:I\rightarrow A^{\sharp}\Box A$ and $e_A:A\Box A^{\sharp}\rightarrow I$, and 2-isomorphisms\footnote{The notation $a^{\bullet}_{A,A^{\sharp},A}$ is used in \cite{SP} to refer to the chosen ajoint pseudo-inverse of the 1-morphism $a_{A,A^{\sharp},A}$.} $$C_A:(e_A\Box A)\circ a^{\bullet}_{A,A^{\sharp},A}\circ (A\Box i_A)\Rightarrow Id_A,$$
$$D_A: Id_{A^{\sharp}}\Rightarrow (A^{\sharp}\Box e_A)\circ a_{A^{\sharp},A,A^{\sharp}}\circ (i_A\Box A^{\sharp}),$$
satisfying the two swallowtail equations (depicted in figures 3 and 4 of \cite{Pstr}). We will also say that the data $(A,A^{\sharp},i_A,e_A,C_A,D_A)$ is a coherent dual pair, or that $(A,A^{\sharp},i_A,e_A,C_A,D_A)$ is a coherent left dual for $A^{\sharp}$.
\end{Definition}

\begin{Remark}
In the notation of definition \ref{def:coherentdual}, if we assume that $\mathfrak{C}$ is a strict cubical monoidal 2-category, the swallowtail equations simplify to $$\left[e_A\circ (C_A\Box A^{\sharp})\right]\cdot \left[\phi_{(e_A,Id_I),(Id_{A\Box A^{\sharp}},e_A)}\circ (A\Box i_A\Box A^{\sharp})\right]\cdot \left[e_A\circ (A\Box D_A)\right] = Id_{e_A},$$
$$\left[ (A^{\sharp}\Box C_A)\circ i_A\right]\cdot \left[(A^{\sharp}\Box e_A\Box A)\circ \phi_{(i_A,Id_I),(Id_{A^{\sharp}\Box A},i_A)}\right]\cdot\left[(D_A\Box A)\circ i_A\right]= Id_{i_A}.$$
\end{Remark}

It is clear that every coherent right dual is a right dual. Hence, it is natural to ask whether every right dual can be made into a coherent right dual. This question was solved in \cite{Pstr}.

\begin{Corollary}{\cite[Corollary 2.8]{Pstr}}
Every right dual can be made coherent, and every left dual can be made coherent.
\end{Corollary}

\begin{Definition}
Let $\mathfrak{C}$ be a monoidal 2-category. We say that $\mathfrak{C}$ has right duals, resp. left duals, if every object has a right dual, resp. left dual. We say that $\mathfrak{C}$ is rigid if it has right and left duals.
\end{Definition}

Using the above result from \cite{Pstr}, we believe that one can construct a right dual 2-functor on any monoidal 2-category with right duals.

\begin{Notation}
Given $\mathfrak{C}$ a monoidal 2-category with monoidal product $\Box$, we denote by $\mathfrak{C}^{\Box op}$ the monoidal 2-category with the opposite monoidal product.
\end{Notation}

\begin{Conjecture}\label{con:2dualfunctor}
Let $\mathfrak{C}$ be a monoidal 2-category with right duals. There exists a monoidal 2-functor $$(-)^{\sharp}:\mathfrak{C}\rightarrow\mathfrak{C}^{\Box op;1op}$$ that sends an object $A$ to the object underlying a right dual for $A$. If, in addition, $\mathfrak{C}$ has left duals, there is a monoidal 2-functor, which we denote by $^{\sharp}(-)$ that sends an object to its left dual. Further, $^{\sharp}(-)$ is a monoidal pseudo-inverse for $(-)^{\sharp}$.
\end{Conjecture}
\begin{Remark}
We prove in lemma \ref{lem:duals2functor} below that the underlying 2-functor exists. However, checking the coherence axioms for a monoidal 2-category involves making sure that big composites of interchangers agree, and we have not found a satisfactory way to deal with these.
\end{Remark}

The decategorified version of the next lemma is well-known.

\begin{Lemma}\label{lem:adjunctionduals}
Let $\mathfrak{C}$ be a monoidal 2-category, and $A$, $B$, and $C$ be objects of $\mathfrak{C}$.
\begin{enumerate}
\item If $C$ has a right dual, there are natural equivalences

$$Hom_{\mathfrak{C}}(A,B\Box C)\simeq Hom_{\mathfrak{C}}(A\Box C^{\sharp},B),$$
$$Hom_{\mathfrak{C}}(C\Box A,B)\simeq Hom_{\mathfrak{C}}(A,C^{\sharp}\Box B).$$

\item If $C$ has a left dual, there are natural equivalences

$$Hom_{\mathfrak{C}}(A\Box C,B)\simeq Hom_{\mathfrak{C}}(A,B\Box ^{\sharp}C),$$
$$Hom_{\mathfrak{C}}(A,C\Box B)\simeq Hom_{\mathfrak{C}}(^{\sharp}C\Box A,B).$$
\end{enumerate}
\end{Lemma}
\begin{proof}
Without loss of generality, we may assume that $\mathfrak{C}$ is strict cubical (this follows from the coherence theorem of \cite{GPS}). Let $(C,C^{\sharp},i_C,e_C,C_C,D_C)$ be a coherent dual pair. Then, the functors

\begin{center}
\begin{tabular}{ c c c }
$Hom_{\mathfrak{C}}(A,B\Box C)$ & $\rightleftarrows$ & $Hom_{\mathfrak{C}}(A\Box C^{\sharp},B)$ \\ 
$f$ & $\mapsto$ & $(B\Box e_C)\circ(f\Box C^{\sharp})$ \\  
$(g\Box C)\circ (A\Box i_C)$ & $\mapsfrom$ & $g$   
\end{tabular}
\end{center}

form an adjoint equivalence with counit $$\lbrack(A\Box C_C)\circ f\rbrack \cdot \lbrack (B\Box e_C\Box C) \circ\phi_{(f,Id),(Id,i_C)}\rbrack$$ and unit $$\lbrack\phi_{(g,Id),(Id,e_C)} \circ (A\Box i_C\Box C^{\sharp})\rbrack \cdot \lbrack g\circ (B\Box D_C)\rbrack.$$ The triangle identites follow from the swallowtail equations. The naturality in $A$ and $B$ is clear from the definition.
\end{proof}

\begin{Remark}
The above lemma can be reformulated by saying that certain 2-functors form a 2-adjunction. For instance, if $C$ has a right dual, then $(-)\Box C^{\sharp}$ is left 2-adjoint to $(-)\Box C$.
\end{Remark}

\subsection{Interactions}

The monoidal 2-categories we will consider have both adjoints and duals. That is why we now examine the properties of such monoidal 2-categories.

\begin{Lemma}\label{lem:dualsandadjoints}
Let $A^{\sharp}$ be a right dual for $A$ in $\mathfrak{C}$, a monoidal 2-category with right adjoints, with unit $i_A$ and counit $e_A$. Then $A^{\sharp}$ has a right dual.
\end{Lemma}
\begin{proof}
Using the monoidal 2-functor $(-)^*$ constructed in lemma \ref{lem:adjointsmonoidal2functor}, one gets that the image of $i_A$ and $e_A$ under $(-)^*$ witness that $A^{\sharp}$ is a right dual for $A$ in $\mathfrak{C}^{1op;2op}$. This means that $A$ is a right dual for $A^{\sharp}$ in $\mathfrak{C}$.
\end{proof}

\begin{Remark}
Lemma \ref{lem:dualsandadjoints} can be generalized to $n$-categories, see \cite[Lem. 4.1.2]{Araujo}.
\end{Remark}

\begin{Corollary}\label{cor:rightdualsrightadjoints}
Let $\mathfrak{C}$ be a monoidal 2-category with right adjoints and right duals. For any object $A$ of $\mathfrak{C}$, $A^{\sharp\sharp}$ is equivalent to $A$.
\end{Corollary}

\begin{Corollary}\label{cor:rightdualsrigid}
Let $\mathfrak{C}$ be a monoidal 2-category with right adjoints. If $\mathfrak{C}$ has right duals, then it also has left duals.
\end{Corollary}

\section{Fusion 2-Categories}\label{sec:fusion}

Throughout, we work over a fixed algebraically closed field $\mathds{k}$ of characteristic zero.

\subsection{(Multi)Fusion 2-Categories}

\begin{Definition}\label{def:f2c}
A multifusion 2-category $\mathfrak{C}$ is a finite semisimple 2-category (\cite{D} definition 2.1.3) equipped with a rigid $\mathds{k}$-linear monoidal structure. In particular, it comes equipped with a bilinear 2-functor $$\Box:\mathfrak{C}\times\mathfrak{C}\rightarrow \mathfrak{C},$$ and a monoidal unit $I$. A fusion 2-category is a multifusion 2-category whose monoidal unit is simple.
\end{Definition}

\begin{Remark}
By their very definition, all the standard results of monoidal 2-category theory (up to linearization) apply to multifusion 2-categories. For instance, by \cite{SP}, every multifusion 2-category is equivalent to a skeletal multifusion 2-category.
\end{Remark}

It is well-known that the monoidal unit of a multifusion category splits as a direct sum of non-isomorphic simple objects (see \cite{EGNO} section 4.3). Using the fact that every object of a finite semisimple 2-categories decomposes into a direct sum of simple objects (see proposition 1.4.5 of \cite{DR} and lemma 2.1.5 of \cite{D}), a similar result holds for multifusion 2-categories.

\begin{Lemma}\label{lem:splitting2multifusion}
Let $\mathfrak{C}$ be a multifusion 2-category. Let $X_i$, $i=1,..n$ be the finitely many simple objects appearing in the decomposition of the monoidal unit as a direct sum of simple objects, i.e. $$I\simeq\boxplus_{i=1}^n X_i.$$ Then, $Hom_{\mathfrak{C}}(X_i,X_j)$ is non-zero if and only if $i=j$. Moreover, we have that $X_i\Box X_j$ is equivalent to $X_i$ if $i=j$ and to $0$ otherwise.
\end{Lemma}
\begin{proof}
For any $i$, we have $X_i\Box I\simeq X_i$. As $X_i$ is simple, there exists precisely one $j$ such that $X_i\Box X_j$ is non-zero. Together with the reverse argument on $j$, this shows that $$X_i \simeq X_i\Box X_j\simeq X_j.$$ If $i\neq j$, then $X_i\Box I$ would have $(X_i\Box X_i) \boxplus (X_i\Box X_j)$ as a summand, whence would not be simple. Thus, we must have $i=j$ and $X_i\Box X_i\simeq X_i$. Moreover, this shows that $X_i$ is both a left and a right dual for $X_i$.

Let $i$, $j$ be arbitrary. Then, we have: $$Hom_{\mathfrak{C}}(X_i,X_j)\simeq Hom_{\mathfrak{C}}(X_i\Box (^{\sharp}X_j), I)\simeq Hom_{\mathfrak{C}}(\delta_{ij}X_i, I).$$ The last term is non-zero precisely when $i=j$. This finishes the proof.
\end{proof}

\begin{Remark}
Lemma \ref{lem:splitting2multifusion} can be seen as a generalization of the fact that braided multifusion categories have no non-zero entries away from the diagonal (see lemma 5.3 of \cite{BDSPV}).
\end{Remark}

Let $\mathfrak{C}$ be a multifusion 2-category. We write $_{i}\mathfrak{C}_{j}$ for the semisimple 2-category $X_i\Box \mathfrak{C}\Box X_j$. The following result is a direct analogue of the usual decomposition of a multifusion category.

\begin{Lemma}
The semisimple 2-categories $_{i}\mathfrak{C}_{i}$ are fusion 2-categories, and the finite semisimple 2-categories $_{i}\mathfrak{C}_{j}$ are $(_{i}\mathfrak{C}_{i}, {_{j}\mathfrak{C}_{j}})$-bimodule 2-categories. Finally, the following matrix
$$\begin{pmatrix}
_{1}\mathfrak{C}_1 & \cdots & _1\mathfrak{C}_n\\
\vdots & \ddots & \vdots\\
_n\mathfrak{C}_1&\cdots & _n\mathfrak{C}_n\\
\end{pmatrix}$$
represents the fusion rule of the multifusion 2-category $\mathfrak{C}$.
\end{Lemma}

\begin{Example}
The 2-category of representations of a finite 2-groupoid $\mathcal{G}$ is semisimple and finite by \cite{DR}. It inherits a (symmetric) monoidal structure from the symmetric monoidal structure on $\mathbf{2Vect}$. The monoidal unit is given by the constant 2-functor $\mathcal{G}\rightarrow \mathbf{2Vect}$ with value $\mathbf{Vect}$. This 2-representation splits as the direct sum of the simple 2-representations that are constant with value $\mathbf{Vect}$ on exactly one component of $\mathcal{G}$ and $0$ on the others.
\end{Example}

\subsection{Comparison with strict Fusion 2-Categories}

Douglas and Reutter have used in \cite{DR} the term fusion 2-category to refer to certain Gray monoids; We call such objects strict fusion 2-categories. Now, we want to compare their definition with ours.

\begin{Lemma}
Let $\mathfrak{C}$ be a multifusion 2-category. There exists a multifusion 2-category $\mathfrak{D}$, whose underlying monoidal 2-category is a strict cubical $\mathds{k}$-linear monoidal 2-category, that is linearly equivalent to $\mathfrak{C}$. Moreover, if $\mathfrak{C}$ is fusion, so is $\mathfrak{D}$.
\end{Lemma}
\begin{proof}
Using a $\mathds{k}$-linear version of the coherence theorem of \cite{GPS}, we obtain a strict cubical $\mathds{k}$-linear monoidal 2-category $\mathfrak{D}$ that is linearly monoidally equivalent to $\mathfrak{C}$. Observe that the underlying linear equivalence of 2-categories witnesses that the 2-category $\mathfrak{D}$ is a finite semisimple 2-category. Moreover, rigidity is preserved by monoidal equivalences of 2-categories. This proves the first part of the result. The last part follows from the Whitehead theorem for monoidal 2-categories (see \cite{SP}).
\end{proof}

\begin{Lemma}\label{lem:strictification1}
There is a bijection between weak fusion 2-categories, whose underlying monoidal 2-category is a strict cubical (or opcubical) $\mathds{k}$-linear monoidal 2-category and strict fusion 2-categories.\footnote{This statement can be made rigorous using a set-theoretic argument. For instance, one could use a bigger universe.}
\end{Lemma}

\begin{proof}
The $\mathds{k}$-linear version of \cite{BMS} lemma 2.16 shows that there is a bijection between strict cubical $\mathds{k}$-linear monoidal 2-categories and $\mathds{k}$-linear Gray monoids with one object. Moreover, the equivalence does not affect the underlying 2-categories, and by lemma 2.1.4 of \cite{D} its finite semisimple in the sense of \cite{DR}. Thus, the only thing we have to prove is that this bijection respects the existence of duals. This property follows from the fact that the monoidal product with a fixed 1-morphism is invariant under this bijection by construction.
\end{proof}

\begin{Remark}
In particular, we may invoke all the results that \cite{DR} have proven for strict fusion 2-categories, and apply them to fusion 2-categories.
\end{Remark}

\subsection{Duals in Multifusion 2-Categories}

Specializing lemma  \ref{lem:dualsandadjoints} to multifusion 2-categories, we obtain the following lemma.

\begin{Lemma}\label{lem:m2cdoubledual}
Let $\mathfrak{C}$ be a multifusion 2-category, and let $A^{\sharp}$ be a right dual for $A$ in $\mathfrak{C}$, then $A$ is a right dual for $A^{\sharp}$.
\end{Lemma}

On the one hand, the decategorified analogue of the next result is well-known: it says that left and right duals in a fusion category agree. The proof relies crucially on the category being semisimple. On the other hand, in the context of fusion 2-categories, the proof has a very distinct flavour; it uses lemma \ref{lem:m2cdoubledual}, which applies in great generality (see corollary \ref{cor:rightdualsrightadjoints}).

\begin{Corollary}
Let $\mathfrak{C}$ be a multifusion 2-category, and $A$ an object of $\mathfrak{C}$. Then, $A^{\sharp\sharp}$ is equivalent to $A$.
\end{Corollary}

\begin{Corollary}\label{cor:monoidalrightdualsss2c}
Let $\mathfrak{C}$ be a monoidal finite semisimple 2-category. If $\mathfrak{C}$ has right duals, then it also has left duals, i.e it is multifusion.
\end{Corollary}

As is the case in any multifusion category, the right dual of a simple object is again a simple object.

\begin{Lemma}\label{lem:simpleduals}
In any multifusion 2-category, the left and right duals of a simple object are simple.
\end{Lemma}
\begin{proof}
Let $A$ be a simple object with right dual $A^{\sharp}$. Observe that the right dual of a non-zero object has to be non-zero. Further, a right dual for a direct sum is given by the direct sum of the right duals. Thus, if $A^{\sharp}$ were not simple, i.e. had two non-zero summands, then $A^{\sharp\sharp}$ would have two non-zero summands. This contradicts the fact that $A\simeq A^{\sharp\sharp}$ is simple.
\end{proof}

Lemma \ref{lem:simpleduals} implies that the operation of taking the right dual induces a bijection on the set of equivalence classes of simple objects. Now, observe that lemma \ref{lem:adjunctionduals} also applies to multifusion 2-categories, yielding the following results:

\begin{Corollary}\label{cor:rightdualcomponent}
Let $\mathfrak{C}$ be a multifusion 2-category, and $A$, $B$ two simple objects such that $Hom_{\mathfrak{C}}(A,B)$ is non-trivial. Then, $Hom_{\mathfrak{C}}(A^{\sharp},B^{\sharp})$ is non trivial.
\end{Corollary}
\begin{proof}
Note that it is enough to prove that $Hom_{\mathfrak{C}}(B^{\sharp},A^{\sharp})$ is non-trivial. Namely, the 2-functor $(-)^*$ of lemma \ref{lem:adjoints2functors} provides us with a linear equivalence:
$$Hom_{\mathfrak{C}}(A^{\sharp},B^{\sharp})\simeq Hom_{\mathfrak{C}}(B^{\sharp},A^{\sharp}).$$
Now, using lemma \ref{lem:adjunctionduals}, there are linear equivalences:
$$Hom_{\mathfrak{C}}(B^{\sharp},A^{\sharp})\simeq Hom_{\mathfrak{C}}(A\Box B^{\sharp},I) \simeq Hom_{\mathfrak{C}}(A,B).$$
This concludes the proof.
\end{proof}

\begin{Corollary}\label{cor:multiplicationcomponent}
Let $\mathfrak{C}$ be a fusion 2-category, and $A$, $J$ two simple objects such that $J$ is in the component of the monoidal unit (i.e. $Hom_{\mathfrak{C}}(I,J)$ is non-zero), then $J\Box A$ is a direct sum of simple objects in the connected component of $A$.
\end{Corollary}
\begin{proof}
Let $B$ be a simple summand of $J\Box A$. There exists a non-trivial 1-morphism between $J\Box A$ and $B$. Thence, there is a non-trivial 1-morphism between $ B\Box A^{\sharp}$ and $J$. By proposition 1.2.19 of \cite{DR}, this implies that there exists a non-trivial 1-morphism between $ B\Box A^{\sharp}$ and $I$, which is equivalent to saying that $A$ and $B$ are in the same connected component.
\end{proof}

We examine the behaviour of simple objects under the monoidal product and arbitrary 2-functors.

\begin{Lemma}\label{lem:fusionproductsimplenonzero}
Let $\mathfrak{C}$ be a fusion 2-category and $C$, $D$ two non-zero objects. Then $C\Box D$ is non-zero.
\end{Lemma}
\begin{proof}
By rigidity, $D$ has a right dual $D^{\sharp}$. In particular, the decomposition of $D\Box D^{\sharp}$ into simple objects contains a copy of $J$, a simple object in the connected component of $I$. By definition, there exists a non-zero 1-morphism $f:I\rightarrow J$. Thus, we get a map $$\xymatrix{C\Box I\ar[r]^-{C\Box f}& C\Box J\ \ar@{^{(}->}[r]& C\Box D\Box D^{\sharp}.}$$ On one hand, if $C\Box D $ were equivalent to zero, then $C\Box J \simeq 0$, whence we would have $C\Box f \simeq 0$. On the other hand, $f$ has a left adjoint $^*f$, and the 2-functor $C\Box (-)$ preserves adjunctions. As $Id_I$ is a direct summand of $^*f\circ f$, we find that $C\Box f\neq 0$. Consequently, $C\Box D$ must be non-zero.
\end{proof}

\begin{Proposition}
Let $F:\mathfrak{C}\rightarrow \mathfrak{D}$ be a monoidal 2-functor from a fusion 2-category to a multifusion 2-category. For any non-zero object $C$ of $\mathfrak{C}$, we have that $F(C)$ is non-zero.
\end{Proposition}
\begin{proof}
As $F$ is a monoidal 2-functor, we know that $F(I)$ is non-zero. Now, the evaluation 1-morphism $e_C:C\Box C^{\sharp}\rightarrow I$ is non-zero, and has a left adjoint. As $F$ preserves adjunctions, we find that $F(e_C)$ is non-zero. Given that $F(C)\Box F(C^{\sharp})\simeq F(C\Box C^{\sharp})$, this implies that $F(C)$ is non-zero.
\end{proof}

\subsection{Connected Fusion 2-Categories}
The goal of this section is to study a special class of fusion 2-categories: connected fusion 2-categories. They are ubiquitous both because the 2-category of finite semisimple module categories associated to a braided fusion category is connected, and because every fusion 2-category has a connected fusion 2-category as a full sub-2-category. This is similar to the fact that every topological monoid admits a connected submonoid given by the connected component of the identity.

\begin{Definition}
A connected fusion 2-category is a fusion 2-category whose underlying finite semisimple 2-category is connected as in definition 1.2.22 of \cite{DR}.
\end{Definition}

\begin{Remark}
By the categorical Schur lemma, i.e. proposition 1.2.19 of \cite{DR}, in order to show that a fusion 2-category is connected, it is enough to check that the $Hom$-categories from the monoidal unit to any simple object is non-trivial.
\end{Remark}

Given a fusion category $\mathcal{C}$, we write $\mathbf{Mod}(\mathcal{C}):=Cau(\mathrm{B}\mathcal{C})$ for the Cauchy completion of $\mathrm{B}\mathcal{C}$ in the sense of \cite{GJF}. By theorem 3.1.7 of \cite{GJF}, we may also think of $\mathbf{Mod}(\mathcal{C})$ as the 2-category of separable algebras, bimodules, and bimodule maps in $\mathcal{C}$. Thus, $\mathbf{Mod}(\mathcal{C})$ is a finite semisimple 2-category by theorem 1.4.8 of \cite{DR}, and it is connected by proposition 2.3.5 of \cite{D}. Further, by proposition 1.3.13 of \cite{DR}, $\mathbf{Mod}(\mathcal{C})$ is equivalent to the 2-category of finite semisimple right $\mathcal{C}$-module categories. The equivalence sends a separable algebra $A$ in $\mathcal{C}$ to $Mod_{\mathcal{C}}(A)$ the finite semisimple category of right $A$-modules in $\mathcal{C}$. Now, if we equip $\mathcal{C}$ with a braiding, more can be said.

\begin{Proposition}{\cite[Construction 2.1.19]{DR}}\label{prop:braidefusionf2c}
Let $\mathcal{C}$ be a braided fusion category. Then, $\mathbf{Mod}(\mathcal{C})$ is a connected fusion 2-category, with monoidal product given by $\boxtimes_{\mathcal{C}}$ the balanced Deligne tensor product.
\end{Proposition}
\begin{proof}
Note that the monoidal 2-category $\mathrm{B}\mathcal{C}$ is rigid. Thence, through the proof of theorem 4.1.1 of \cite{GJF}, we find that its Cauchy completion, $\mathbf{Mod}(\mathcal{C})$, is a multifusion 2-category with monoidal product $\Box$. Explicitly, the monoidal structure is as follows: Given two separable algebras $A$, $B$ in $\mathcal{C}$, representing two right $\mathcal{C}$-module categories $Mod_{\mathcal{C}}(A)$ and $Mod_{\mathcal{C}}(B)$, their product $A\otimes B$ is again a separable algebra in $\mathcal{C}$. The separable algebra $A\otimes B$ represents an object $Mod_{\mathcal{C}}(A\otimes B)$ in $\mathbf{Mod}(\mathcal{C})$, which is, by construction, the monoidal product of $Mod_{\mathcal{C}}(A)$ and $Mod_{\mathcal{C}}(B)$. The result follows from the equivalence of right $\mathcal{C}$-module categories:
$$Mod_{\mathcal{C}}(A)\boxtimes_{\mathcal{C}} Mod_{\mathcal{C}}(B)\simeq Mod_{\mathcal{C}}(A)\boxtimes_{\mathcal{C}} RMod_{\mathcal{C}}(B^{op})$$
$$\simeq Bimod_{\mathcal{C}}(A, B^{op})\simeq Mod_{\mathcal{C}}(A\otimes B) = Mod_{\mathcal{C}}(A)\Box Mod_{\mathcal{C}}(B),$$
where we have used the equivalence of right $\mathcal{C}$-module categories $$Mod_{\mathcal{C}}(B)\simeq RMod_{\mathcal{C}}(B^{op})$$
between the category of left $B$-modules and the category of right $B^{op}$-modules.
\end{proof}

\begin{Definition}
Let $\mathfrak{C}$ be a fusion 2-category with monoidal unit $I$. We denote by $\mathfrak{C}^0$ the connected component of the identity, i.e. the full additive sub-2-category on the simple objects that admit a non-zero morphism from $I$.
\end{Definition}

\begin{Proposition}
The 2-category $\mathfrak{C}^0$ is a fusion sub-2-category of $\mathfrak{C}$ that is connected.
\end{Proposition}
\begin{proof}
We begin by proving that $\mathfrak{C}^0$ is finite semisimple. The only property which is not obvious is that $\mathfrak{C}^0$ has all condensates. Let $(A,...)$ be a 2-condensation monad in $\mathfrak{C}^0$, and let $(A,B,f,g,...)$ be an extension in $\mathfrak{C}$ of $(A,...)$ to a 2-condensation. Observe that for every simple summand $C$ in $\mathfrak{C}$ of $B$, the composite of $f$ with the projection $B\rightarrow C$ is a non-zero 1-morphism $A\rightarrow C$. (If this 1-morphism was zero, $(A,B,f,g,...)$ would not be a 2-condensation.) Thus, we find that $C$ is in $\mathfrak{C}^0$. Further, by definition, we have that $\mathfrak{C}^0$ is connected.

By corollary \ref{cor:multiplicationcomponent}, the monoidal product of $\mathfrak{C}$ restricts to give $\mathfrak{C}^0$ a monoidal structure. Finally, as $I^{\sharp}\simeq I$, we find by corollary \ref{cor:rightdualcomponent} that this monoidal structure is rigid, and $\mathfrak{C}^0$ is clearly fusion.
\end{proof}

\begin{Corollary}\label{cor:connectedcomponent}
Let $\mathfrak{C}$ be a fusion 2-category. Then, there is an equivalence of monoidal 2-categories:
$$Mod(End_{\mathfrak{C}}(I))\simeq \mathfrak{C}^0.$$
\end{Corollary}
\begin{proof}
As a consequence of the proof of the above proposition, we find that $$\mathrm{B}End_{\mathfrak{C}}(I)\hookrightarrow \mathfrak{C}^0$$ is a Cauchy completion. Further, this inclusion is monoidal, whence, by the 3-universal property of the Cauchy completion, we get the desired result.
\end{proof}

Corollary \ref{cor:connectedcomponent} shows that the behavior of the monoidal product on the connected component of the identity is completely determined by the braiding on the fusion category $End_{\mathfrak{C}}(I)$.

\begin{Proposition}\label{prop:braidedfusioncatconnectedf2c}
There is an equivalence between the category of connected fusion 2-categories and equivalence classes of monoidal linear 2-functors, and the category of braided fusion categories and equivalence classes of braided tensor functors.
\end{Proposition}
\begin{proof}
Let us denote by $\mathcal{A}$ the category of connected fusion 2-categories and equivalence classes of monoidal linear 2-functors, and by $\mathcal{B}$ the category of braided fusion categories and equivalence classes of braided tensor functors. Taking the endomorphism category of the monoidal unit yields a functor $$End_{(-)}(I):\mathcal{A}\rightarrow \mathcal{B},$$ and taking the Cauchy completion of the delooped braided fusion category gives a functor $$\mathbf{Mod}(-)=Cau(\mathrm{B}(-)):\mathcal{B}\rightarrow \mathcal{A}.$$ Using the 3-universal property of the Cauchy completion, one finds that these functors are pseudo-inverses for one another.
\end{proof}

\subsection{Examples given by Pointed Braided Fusion Categories}

Let $\mathcal{C}$ be a pointed braided fusion category. By results of $\cite{EGNO}$, we know that this corresponds equivalently to the data of a finite abelian group $A$ equipped with an abelian 3-cocycle $(\omega,\beta)$. We will assume that $\omega$ is trivial. (If $A$ has odd order, this can always be done.) We denote the braided fusion category associated to this data by $\mathbf{Vect}_A^{\beta}$. It is known that finite semisimple indecomposable right module categories over $\mathbf{Vect}_A^{\beta}$ correspond to pairs $(E,\phi)$, where $E$ is a subgroup of $A$ and $\phi$ is 2-cocycle on $E$ with value in $\mathds{k}^{\times}$ (considered up to 2-coboundary). We denote the corresponding right module category by $\mathcal{M}(E,\phi)$. Proposition 3.16 of \cite{ENO2}, explains how to compute the relative Deligne tensor product of two right $\mathcal{C}$-modules when $\beta$ is trivial. We now generalize this result.

Let $(E,\phi)$, and $(F,\psi)$ be two pairs consisting of a subgroup of $A$, and an appropriate 2-cocycle. Let $Alt(\phi):E\times E\rightarrow \mathds{k}^{\times}$ and $Alt(\psi):F\times F\rightarrow \mathds{k}^{\times}$ be the corresponding skew-symmetric bilinear forms, i.e. $$Alt(\phi)(e_1,e_2):= \phi(e_1,e_2)/\phi(e_2,e_1),$$ and similarly for $Alt(\psi)$. We define a skew-symmetric bicharacter $b$ on $E\oplus F$ by $$b((e_1,f_1),(e_2,f_2)):= Alt(\phi)(e_1,e_2)Alt(\psi)(f_1,f_2)\beta(f_1,e_2)/\beta(f_2,e_1).$$ The group $E\cap F$ embeds in $E\oplus F$ via $e\mapsto (e,-e)$, thus we can consider its orthogonal complement $(E\cap F)^{\bot}$ under the bicharacter $b$. Now, let $H$ be the image of $(E\cap F)^{\bot}$ under the canonical map $E\oplus F\rightarrow A$. The restriction of $b$ to $(E\cap F)^{\bot}$ descends to a skew-symmetric bilinear form $b'$ on $H$, which corresponds to an element of $H^2(H,\mathds{k}^{\times})$ represented by a chosen 2-cocycle $\rho$ (see proposition 3.6 of \cite{Tam}).

\begin{Proposition}\label{prop:fusionrulepointedbraided}
Write $\mathcal{C}=\mathbf{Vect}_A^{\beta}$, we have
$$\mathcal{M}(E,\phi)\boxtimes_{\mathcal{C}}\mathcal{M}(F,\psi) \simeq \boxplus_{i=1}^m \mathcal{M}(H,\rho),$$
where $$m=\frac{|(E\cap F)^{\bot}||(E\cap F)|}{|E|| F|}.$$
\end{Proposition}
\begin{proof}
Let $A(E,\phi)$ be the algebra in $\mathcal{C}$ whose underlying object is $\mathds{k}E$ and whose multiplication is given by $\phi$, and similarly for $A(F,\psi)$. By definition, $\mathcal{M}(E,\phi) = Mod_{\mathcal{C}}(A(E,\phi))$, and $\mathcal{M}(F,\psi) = Mod_{\mathcal{C}}(A(F,\psi))$. The relative Deligne tensor product is given by the category of left modules over the algebra $A(E,\phi)\otimes A(F,\psi)$ in $\mathcal{C}$. Note that the multiplication of this algebra is twisted by the braiding $\beta$ of $\mathcal{C}$. More precisely, the multiplication is given by the 2-cocycle $\tau$ on $E\oplus F$ defined by $$((e_1,f_1),(e_2,f_2))\mapsto  \phi(e_1,e_2)\psi(f_1,f_2)\beta(f_1,e_2).$$ In particular, the corresponding skew-symmetric bicharacter on $E\oplus F$ is $b$ as defined above. Finally, using proposition 2.11 of \cite{ENO2}, we obtain the desired result.
\end{proof}

\begin{Remark}
It should be possible to generalize proposition \ref{prop:fusionrulepointedbraided} to the case were $\omega$ is not assumed to be trivial. However, the above proof does not immediately generalize because if $\omega$ is not trivial, then $\phi$ and $\psi$ may not be 2-cocycles, and thus $Alt(\phi)$ and $Alt(\psi)$ may not be skew-symmetric bicharacters.
\end{Remark}

\begin{Example}\label{ex:Mod(p)}
Recall the notation of example 2.1.12 of \cite{D}. It is well-known that the fusion category $\mathbf{Vect}_{\mathbb{Z}/p\mathbb{Z}}$ admits $p$ distinct braided structures (up to braided monoidal automorphisms of $\mathbf{Vect}_{\mathbb{Z}/p\mathbb{Z}}$ that are the identity on objects). A braiding $b$ on $\mathbf{Vect}_{\mathbb{Z}/p\mathbb{Z}}$ is determined by its value $b_{\mathds{k}_1, \mathds{k}_1} = e^{\frac{2\pi i k}{p}}$, for any $0\leq k<p$. If we allow arbitrary braided monoidal automorphisms, there are two distinct braided structures if $p=2$, and three otherwise. The symmetric or trivial braiding is specified by $k=0$, and we denote the corresponding braided fusion category by $\mathbf{Vect}_{\mathbb{Z}/p\mathbb{Z}}^{triv}$. The other ones correspond to the cases where $0<k<p$ is either a quadratic residue or not. For simplicity, we only treat the case where $k$ is a quadratic residue, for which we may assume $k=1$, and denote the corresponding braided fusion category by $\mathbf{Vect}_{\mathbb{Z}/p\mathbb{Z}}^{\beta}$. The other case is entirely analogous.

Let us begin by examining the monoidal product on the finite semisimple 2-category of finite semisimple right modules categories over $\mathcal{C}:=\mathbf{Vect}_{\mathbb{Z}/p\mathbb{Z}}^{triv}$. As $\mathcal{C}$ is the monoidal unit for the induced monoidal structure on $\mathbf{Mod}(\mathcal{C})$, we only have to determine $\mathbf{Vect}\boxtimes_{\mathcal{C}}\mathbf{Vect}$. Recall that as right $\mathbf{Vect}_{\mathbb{Z}/p\mathbb{Z}}$-module categories, we have $\mathbf{Vect}\simeq \mathcal{M}(\mathbb{Z}/p\mathbb{Z},triv)$ in the notations used above. Now, a straightforward computation using proposition \ref{prop:fusionrulepointedbraided} shows: $$\mathbf{Vect}\boxtimes_{\mathcal{C}}\mathbf{Vect}\simeq \boxplus_{i=1}^p\mathbf{Vect},$$
as right $\mathcal{C}$-module categories.

We now turn our attention to the case $\mathcal{D}:=\mathbf{Vect}_{\mathbb{Z}/p\mathbb{Z}}^{\beta}$. As above, $\mathcal{D}$ is the monoidal unit of the induced monoidal structure on $\mathbf{Mod}(\mathcal{D})$, whence we only have to determine $\mathbf{Vect}\boxtimes_{\mathcal{D}}\mathbf{Vect}$. In order to use proposition \ref{prop:fusionrulepointedbraided}, we compute that $\langle (1,-1)\rangle^{\bot}\subseteq \mathbb{Z}/p\mathbb{Z}\oplus \mathbb{Z}/p\mathbb{Z}$ is precisely $\langle (1,-1)\rangle$. This gives $m=1$, and thus $$\mathbf{Vect}\boxtimes_{\mathcal{D}}\mathbf{Vect}\simeq \mathcal{M}(0,triv) \simeq\mathcal{D}.$$ This examples shows that the braiding we put on a fusion category can have a big impact on the fusion rule of the associated fusion 2-category.
\end{Example}

\appendix

\section{Two 2-Functors}\label{sec:two2functors}

\subsection{The Adjoints Monoidal 2-Functor}\label{sub:adjointsmonoidal2functor}

\begin{Notation}
Let $\mathfrak{C}$ be a (monoidal) 2-category. Given $f:A\rightarrow B$ a 1-morphism in $\mathfrak{C}$, we denote by $\{f\}^{1op;2op}: B\rightarrow A$ the corresponding 1-morphism in $\mathfrak{C}^{1op;2op}$. We write $\circ^{op}$ for the composition of 1-morphisms in $\mathfrak{C}^{1op;2op}$. Given $\alpha:f\Rightarrow g$ a 2-morphism in $\mathfrak{C}$, we denote by $\{\alpha\}^{1op;2op}: \{g\}^{1op;2op}\Rightarrow \{f\}^{1op;2op}$ the corresponding 2-morphism in $\mathfrak{C}^{1op;2op}$.
\end{Notation}

\renewcommand*{\proofname}{Proof of lemma \ref{lem:adjoints2functors}}

\begin{proof}
This is well-known (for instance, see \cite{Ghosh}), but let us indicate briefly how to proceed. Thanks to the coherence theorem for 2-categories, we can omit the coherence 2-isomorphisms for $\mathfrak{C}$. The 2-functor $(-)^*$ is defined as follows:

It sends the obejct $C$ of $\mathfrak{C}$ to itself. Given a 1-morphism $f$, we set $f^*:= \{f'\}^{1op;2op}$, where $f'$ is a fixed right adjoint for $f$ with unit $\eta_f$ and counit $\epsilon_f$. Given a 2-morphism $\alpha:f\Rightarrow g$, we define $$\alpha^*:=\{(f'\circ \epsilon_g)\cdot (f'\circ\alpha\circ g')\cdot (\eta_f\circ g')\}^{1op;2op}.$$ Given two 1-morphisms $f:A\rightarrow B, g:B\rightarrow C$, the structure 2-isomorphism $g^*\circ^{op} f^*\Rightarrow (g\circ f)^*$ in $\mathfrak{C}^{1op;2op}$ witnessing that $(-)^*$ respects the composition of 1-morphisms is given by the 2-isomorphisms $$\{(f'\circ g'\circ \epsilon_{g\circ f})\cdot (f' \circ \eta_g \circ f\circ (g\circ f)')\cdot (\eta_f \circ (f\circ g)')\}^{1op;2op}.$$  Similarly, the unitors are given by $\{\eta_{Id_C}^{-1}\}^{1op;2op}:Id_C^*\Rightarrow Id_C$ for every object $C$ of $\mathfrak{C}$. This data specifies the 2-functor $(-)^*$. The associativity of composition follows from the triangle identities, and the unitality of composition is clear.

Further, if $\mathfrak{C}$ also has left-adjoints, we can dually define a 2-functor $^*(-)$. In particular, for any given 1-morphism $f$, we fix a left adjoint $'f$ with unit $\xi_f$ and counit $\kappa_f$. It is not hard to see that $^*(-)$ is a pseudo-inverse for $(-)^*$ (see the footnote on page \pageref{lem:adjoints2functors}). For instance, a 2-natural equivalence $\theta:(^*(-))^*\Rightarrow Id$ is given by the identity 1-morphism on objects, and on a 1-morphism $f$ by the 2-isomorphism $( \epsilon_{('f)}\circ f )\cdot ( ('f)'\circ\xi_f ).$
\end{proof}

\renewcommand*{\proofname}{Proof}

\begin{Lemma}\label{lem:technicaladjoint2functor}
Let $F:\mathfrak{C}\rightarrow \mathfrak{D}$ be a 2-functor. There is a 2-natural equivalence $e$ that fits into the following diagram:
$$\xymatrix{\mathfrak{C}\ar[r]^-F\ar[d]_-{(-)^*}& \mathfrak{D}\ar[d]^-{(-)^*}\\\mathfrak{C}^{1op;2op}\ar@{=>}[ru]^{e}\ar[r]_-{F^{1op;2op}}&\mathfrak{D}^{1op;2op}.}$$
\end{Lemma}
\begin{proof}
As above, we omit the coherence 2-isomorphisms of $\mathfrak{C}$, and $\mathfrak{D}$. Given an object $C$ of $\mathfrak{C}$, we set $e_C := Id_{F(C)}$. On the 1-morphism $f:C\rightarrow D$, we define $e_f$ as the following 2-isomorphism:

$$\{(\eta_{F(f)}\circ F(f'))\cdot (F(f)'\circ F_{f',f}^{-1}) \cdot (F(f)'\circ F(\epsilon_{f}))\}^{1op;2op},$$

where $F_{f',f}:F(f')\circ F(f)\Rightarrow F(f'\circ f)$ is the coherence 2-isomorphism supplied by $F$. It is not hard to check that this defines a 2-natural equivalence.
\end{proof}

\renewcommand*{\proofname}{Proof of lemma \ref{lem:adjointsmonoidal2functor}}

\begin{proof}
Thanks to the coherence theorem for monoidal 2-cate\-gories, there is an equivalence of monoidal 2-categories $F:\mathfrak{C}\rightarrow\mathfrak{D}$ such that $\mathfrak{D}$ is strict cubical. Below, we will endow the adjoint 2-functor $(-)^*$ of $\mathfrak{D}$ with a monoidal structure. Through the natural 2-equivalence of lemma \ref{lem:technicaladjoint2functor}, this shows that the adjoint 2-functor $(-)^*$ of $\mathfrak{C}$ is equivalent as a 2-functor to a monoidal one. Hence, it is monoidal itself.

In order to specify a monoidal structure on the 2-functor $(-)^*$ on $\mathfrak{D}$, we need to give some data. To make this more digestible, we use the notations of \cite{SP}. We begin by defining the 2-natural equivalence $$\chi:(-)^*\Box (-)^*\Rightarrow ((-)\Box (-))^*.$$ Given two objects $A$, $B$ we let $\chi_{A,B}:=Id_{A\Box B}$. Given two 1-morphisms $f:A\rightarrow B$ and $g:C\rightarrow D$ we let the 2-isomorphism $\chi_{f,g}:f^*\Box g^*\Rightarrow (f\Box g)^*$ be given by

$$\{\big((f'\Box g')\circ \epsilon_{f\Box g}\big)\cdot \big(\phi^{-1}_{(f',g'),(f,g)}\circ (f\Box g)'\big)\cdot \big((\eta_f\Box\eta_g)\circ(f\Box g)'\big)\}^{1op;2op}.$$

It is not hard to see that $\chi$ is a 2-natural transformation. Further, it is clearly an isomorphism; and we pick $\chi^{\bullet}$ to be its inverse.

We choose the 1-equivalence $\iota$ to be the identity 1-morphism on the monoidal unit. The modifications $\omega$, $\gamma$, $\delta$ are uniquely specified by the universal property of right adjoints. Namely, we let $\gamma_C$ be the 2-isomorphism $\{\eta_{Id_C}\}^{1op;2op}$ in $\mathfrak{C}$ and $\delta_C$ be the identity 2-morphism of $\{Id_C\}^{1op;2op}$, for every object $C$ of $\mathfrak{C}$. Further, for every $A,B,C$ in $\mathfrak{C}$, the modification $\omega_{A,B,C}$ is given by the 2-isomorphism $\{\epsilon_{Id_{A\Box B\Box C}}\}^{1op;2op}$.

The commutativity of the coherence diagrams can be checked using the uniqueness up to unique isomorphism of right adjoints. Similarly, one can endow the 2-functor $^*(-)$ with a monoidal structure.

Finally, we need to construct monoidal 2-natural equivalences witnessing that $(-)^*$ and $^*(-)$ are pseudo-inverse monoidal 2-functors. We construct the monoidal 2-natural equivalence $(^*(-))^*\Rightarrow Id$, the other one can be constructed analogously. Observe that using the argument at the beginning of this proof, it is enough to construct this monoidal 2-natural equivalence on $\mathfrak{D}$. As its underlying 2-natural equivalence we take the 2-natural equivalence $\theta$ defined in the proof of lemma \ref{lem:adjoints2functors} above. Then, in the notation of \cite{SP}, the 2-isomorphism $M$ is given by the identity 2-morphism on $Id_I$, and the modification $\Pi$ is specified on $A,B$ in $\mathfrak{C}$ by the 2-isomorphism $\{\eta_{Id_{A\Box B}}\}^{1op;2op}$.
\end{proof}

\renewcommand*{\proofname}{Proof}

\subsection{The Duals 2-Functor}

\begin{Notation}
Let $\mathfrak{C}$ be a 2-category. Given $f:A\rightarrow B$ a 1-morphism in $\mathfrak{C}$, we denote by $\{f\}^{1op}: B\rightarrow A$ the corresponding 1-morphism in $\mathfrak{C}^{1op}$. We write $\circ^{op}$ for the composition of 1-morphisms in $\mathfrak{C}^{1op}$. Given $\alpha:f\Rightarrow g$ a 2-morphism in $\mathfrak{C}$, we denote by $\{\alpha\}^{1op}: \{f\}^{1op}\Rightarrow \{g\}^{1op}$ the corresponding 2-morphism in $\mathfrak{C}^{1op}$.
\end{Notation}

\begin{Lemma}\label{lem:duals2functor}
Let $\mathfrak{C}$ be a monoidal 2-category that has right duals. Then, there is a 2-functor $$(-)^{\sharp}:\mathfrak{C}\rightarrow \mathfrak{C}^{1op}$$
that sends an object $A$ in $\mathfrak{C}$ to the object underlying a right dual for $A$.
\end{Lemma}

\begin{proof}
We may assume that the monoidal 2-category $\mathfrak{C}$ is strict cubical. For every object $A$ in $\mathfrak{C}$, choose a coherent right dual $(A,A^{\sharp}, i_A, e_A, C_A, D_A)$. These choices allow us to define the following assignments towards defining the 2-functor $(-)^{\sharp}$. An object $A$ in $\mathfrak{C}$ is sent to $A^{\sharp}$. A 1-morphism $f:A\rightarrow B$ is sent to $f^{\sharp}:=\{f'\}^{1op}$, where $$f':= (A^{\sharp}\Box e_B) \circ (A^{\sharp}\Box f\Box B^{\sharp}) \circ (i_A\Box B^{\sharp}): B^{\sharp}\rightarrow A^{\sharp}.$$

A 2-morphism $\alpha: f\Rightarrow g:A\rightarrow B$ is sent to $\alpha^{\sharp}:=\{\alpha'\}^{1op}$, where $$\alpha':= (A^{\sharp}\Box e_B) \circ (A^{\sharp}\Box \alpha\Box B^{\sharp}) \circ (i_A\Box B^{\sharp}):f'\Rightarrow g'.$$

Given two 1-morphisms $f:A\rightarrow B$, $g:B\rightarrow C$ in $\mathfrak{C}$, then $\{-\}^{1op}$ of the following 2-isomorphism in $\mathfrak{C}$: $$
\adjustbox{scale=0.95}{$\begin{tikzcd}
f'\circ g' \arrow[Rightarrow, d, "\phi"]                                                                                                                                                                                                                   \\
\Big(A^{\sharp} \Box \big(e_C \circ (g\Box C^{\sharp})\big)\Big)\circ \Big(A^{\sharp}\Box \big((e_B\Box B)\circ (A\Box i_B)\big) \Box C^{\sharp}\Big)\circ \Big(\big( (A^{\sharp}\Box f) \circ i_A\big) \Box C^{\sharp}\Big) \arrow[Rightarrow, d, "D_B^{-1}"] \\
(g\circ f)'                                                                                                                                                                                                                           
\end{tikzcd}$}
$$

serves as the structure 2-isomorphism in $\mathfrak{C}^{1op}$ witnessing that $(-)^{\sharp}$ respects the composition of 1-morphisms. The unitor on $A$ is provided by $\{D_A\}^{1op}$. Using naturality of the interchangers, it is not hard to check the coherence axioms, and so defines a 2-functor.
\end{proof}

\let\clearpage\newpage

\bibliography{bibliography.bib}

\end{document}